\documentclass[12pt]{amsart}
\usepackage[margin=1.0in]{geometry}
\usepackage{mathptmx}
\usepackage[foot]{amsaddr}
\begin{document}
\title{Mathematical Symbolism in a Russian Literary Masterpiece}
\author{Noah Giansiracusa$^*$} 
\address{$^*$Assistant Professor of Mathematics, Swarthmore College} 
\author{Anastasia Vasilyeva$^{**}$} 
\address{$^{**}$Undergraduate Economics Major, Swarthmore College}
\email{ngiansi1@swarthmore.edu, avasily1@swarthmore.edu}

\maketitle

\begin{abstract}
Andrei Bely's modernist novel \emph{Petersburg}, first published in 1913, is considered a pinnacle of the Symbolist movement.  Nabokov famously ranked it as one of the four greatest masterpieces of 20th-century prose.  The author's father, Bugaev, was an influential mathematician and for 12 years served as the president of the Moscow Mathematical Society; he was also a source of inspiration for one of the main characters in the son's novel.  While the philosophical views and political leanings of the mathematicians surrounding Bely, and their impact on \emph{Petersburg}, have been a topic of recent academic interest, there has not yet been a direct investigation of the surprisingly frequent and sophisticated mathematical passages in the book itself.  We attempt here to rectify this gap in the scholarly literature, and in doing so find a rich tapestry of mathematical ideas and allusions.
\end{abstract}

\section{Introduction}

Andrei Bely (1880-1934) was a fascinating, if somewhat tragic, figure.  He grew up under the stern tutelage of his internationally renowned math professor father, Nikolai Bugaev (1837-1903), and studied the natural sciences at the University of Moscow before turning to the literary arts.  In fact, it was because of the father's fame and reputation that Bely took up a pseudonym when launching his literary career, rather than publishing under his birth name, Boris Bugaev \cite[p. 31]{Mochulsky}.  Both father and son were polymaths with vast philosophical interests.  They seemed intent on expanding the theoretical foundations of their respective subjects (math for the father, poetry and prose for the son) to quixotic extents: they believed they could explain, and even justify, historical events and contemporary social movements through their abstract work.  They lived during tumultuous times in Russia, let us not forget.

Bely's personal philosophy, which indelibly imprinted itself upon his creative work, underwent many stages of evolution, multiple times retreating from what he later admitted were wrong turns, and the critical reaction to his poems and novels was continually mixed during his lifetime.  (The professional jealousies and misguided romantic relationships of his life did not help his cause.)  However, over the years since his death, one work of Bely in particular has distinguished itself not just as Bely's greatest creation, but as one of the true landmarks in modern literature.  Indeed, no less a consummate authority than Vladimir Nabokov (1899-1977) declared Bely's novel \emph{Petersburg} one of the four greatest works of prose of the 20th century, along with Joyce's \emph{Ulysses}, Kafka's \emph{Metamorphosis}, and Proust's \emph{\`A la recherche du temps perdu}.  The scholarly literature on \emph{Petersburg} is vast, exploring the novel from a myriad of perspectives, and it continues to thrive (see, for instance, the recent volume of essays celebrating the novel's centennial \cite{Centennial}).

That mathematics plays a role in \emph{Petersburg} is immediately apparent: how many novels can you think of with a section in the first chapter titled ``Squares, Cubes, and Parallelepipeds''?  The two main protagonists are father and son, and they bear many resemblances to the real-life father and son, Bugaev and Bely.  For instance, the fictionalized father, a conservative senator, has a humorous penchant for contemplating geometric figures in order to soothe his nerves.  So what can be learned about Bely's complex masterpiece by turning an eye to mathematics?  We are not the first to ask this question.

In 2009, historian of Russian science Loren Graham released an engaging and thought-provoking book written collaboratively with mathematician Jean-Michel Kantor in which they explore the religious and mystical turns that mathematics took in Russia during the early 20th century \cite{Naming}.  The key mathematical topic is set theory and the study of infinity, a perplexing realm unleashed by Georg Cantor (1845-1918), who was born, aptly enough, in Petersburg, Russia.  The key mathematicians in their book are the students of Bugaev, and the students of these students.  Bely does not feature in their book, except in passing as a classmate of the protagonists, Pavel Florensky (1882-1937) and Nikolai Luzin (1883-1950), and as a generally relevant poet at the time.  However, implicit in Graham-Kantor's writing is the significant insight that some of the same ontological issues these mathematicians were grappling with in their mathematical work must also have played a role in Bely's literary explorations. 

In 2013, professor of comparative literature and linguistics Ilona Svetlikova published a book exploring the metamathematical world of Bugaev and the Moscow Mathematical School that he is credited with creating.  In the introduction she asserts:  ``The link between his [Bugaev's] profession and the senator's geometrical meditations must have seemed too obvious to deserve attention, and the meaning of this link was thus never explored.'' \cite[p. 2]{Svetlikova}.  In this fascinating scholarly book, Svetlikova thoroughly unwinds the (sometimes bizarre) political ideas of the Moscow Mathematical School and, in broad strokes, connects it to Bely's oeuvre.  The implications for \emph{Petersburg} are based more on general themes in the book than the many specific mathematical passages it contains.  This expands upon earlier work of Lena Szil\'ard \cite{Szilard1,Szilard2}, arguably the first scholar to draw a direct connection between Bugaev's philosophy and Bely's writing.  

Elizabeth Kosakowska's Columbia Ph.D. Dissertation \cite{Columbia}, published the same year as Svetlikova's book, turns the focus squarely on Bely's \emph{Petersburg} by providing a literary exegesis through the lens of popular scientific revolutions that were in the air at the time: quantum mechanics, relativity, etc.  Kosakowska includes the discovery of non-Euclidean geometry as one of the influences on Bely, so pure mathematics enters the scene.  However, while her study spends the majority of its pages looking directly in the text of \emph{Petersburg}, it does so primarily to construe the novel in terms of the philosophy of science and mathematics outside of the novel.

What appears to be missing in this recent spate of mathematically-oriented scholarship on Bely is a direct interpretation and contextualization of the mathematics within \emph{Petersburg} itself---and this is precisely what we strive to provide with this article.  Our primary targeted audience is mathematicians and mathematics enthusiasts, as we hope such readers will be excited to learn of remarkable interactions between their discipline and the literary world.  However, aside from a couple places where we construe Bely's imaginative poetic passages in technical mathematical language, we endeavor to be accessible to non-mathematicians as well.  Although we avoid any attempts at explaining the ``meaning'' of the novel, perhaps literary scholars will find useful our historical contextualization of the math in \emph{Petersburg} and the relations we find in it to other mathematically-oriented works of fictions.  Also, similar to Ahearn's discussion in the \emph{American Mathematical Monthly} that Tolstoy's use of calculus in \emph{War and Peace} provides math teachers with a rich opportunity for thought experiments and a non-traditional learning environment \cite{Tolstoy}, we hope that---with the help of this article---the ideas Bely plays with in his novel will find their way into math classes and particularly appeal to the liberal arts minded student.\\

\emph{Acknowledgements}.  We would like to thank Emily Frey\footnote{Visiting Assistant Professor of Russian Languages and Literature, Swarthmore College, and wife of the first author.} for helpful discussions regarding this paper.  The first author was supported in part by an NSA Young Investigator Grant and wishes to thank that organization for financial support.

\section{Andrei Bely, his father, and his novel}

There are many excellent biographical works written about Bely, in addition to his memoirs and semi-autobiographical novels and poems; so let us refer the interested reader to these for generalities on Bely's life and begin instead with a few specific anecdotes.  We present these anecdotes to illustrate that despite Bely's choice to avoid the career path of his father, he nonetheless inherited some amusing characteristics and propensities associated with being a mathematician.  

A visitor to Bely's apartment once remarked \cite[p. 177]{Mochulsky}: ``The room was covered with books and manuscripts.  For some reason, there was a blackboard in it.''  A blackboard!  Arguably, there is no more canonical symbol of our profession than the omnipresence of this object.  The visitor said Bely sketched spirals on the blackboard then declared: ``Do you see?  The lowest point of the spiral?  That is you and I at present.  That is the present moment of the revolution.  We will not fall further.  The spiral goes toward the top, widening.  It is already carrying us out of hell into the open air.''  Bely, the true syncretist, absorbed ideas from all disciplines and combined them in his own unique way.  The spiral was an important symbol for Bely, which he noted is a ``circular line,'' a figure combining linear and circular modes of thought \cite[p. 48]{Columbia}.  Bely was excited by the idea that the spiral breaks an otherwise stagnant circularity of motion and lifts one to higher and higher levels of being.  This adds an appealing philosophical dimension to the picture we all see in undergraduate topology of the real number line sitting over the unit circle as a universal covering space.  For the literary relevance of the spiral to \emph{Petersburg}, see \cite{Columbia}.

Another instance of Bely merging mathematical thinking with literary concerns is in his memoirs, where he describes one person's attitude toward him as ``concave'' since it drew Bely out into increasing development of his ideas, leading him in the direction in which his thought pointed, and contrasted this with the ``convex'' Bugaev who sought instead to impress his own views upon his son \cite[p.8]{Elsworth}.  Bely's relationship with his father, and thus with mathematics, was extremely complicated.  A rather intriguing excerpt on Bely's childhood from his semi-autobiographical novel \emph{Kotik Letaev} is poignantly summarized by Mochulsky as follows:
\begin{quote}
He is five years old; new life experiences take place; he discovers that his development is ``premature and abnormal.'' His mother kisses him and suddenly pushes him away and begins to cry: ``He's not like me: he's like --- his father...'' He is also crying... ``Is it really my fault that I --- know: --- my papa is in correspondence with Darboux; Poincar\'e likes him; but Weierstrasse not particularly; Idealov was in Leipzig with ... an elliptical function.'' \cite[p. 170]{Mochulsky}
\end{quote}
We also find reference to Henri Poincar\'e (1854-1912) in the extensive bibliography appearing in Bely's theoretical work on Symbolism \cite[p. 34]{Svetlikova}.  Possessing immense libraries covering diverse topics was undeniably a commonality between Bely and his father.

Bugaev was not only the father of Bely, he was the father of the Moscow Mathematics School.  His most successful student was Dmitri Egorov (1869-1931), whose mathematical family tree is depicted on \cite[p. 163]{Naming} and includes many names you will surely recognize from a modern mathematics education: Alexandrov, Arnold, Gelfand, Kolmogorov, Lyapunov, Novikov, Pontryagin, Postnikov, Uryson, etc.  Bugaev insisted that the journal of the Moscow Mathematical Society be published in the Russian language \cite{Demidov-150}, a decision that no doubt strongly impacted the development of ``Russian'' mathematics throughout the 20th century.

Bugaev and Poincar\'e both gave talks at the first ever ICM, which took place in Z\"urich in 1897.  This was 39 years before the first Fields Medal was awarded at the ICM, and three years before the seminal Paris ICM talk of David Hilbert (1862-1943) in which he introduced his legendary list of open problems.  Poincar\'e's ICM talk extolled the beauty of pure mathematics and the remarkable success of its application, particularly based on analysis, to physics and celestial mechanics \cite{Poincare-ICM}.  Bugaev's ICM talk \cite{Bugaev-ICM}, on the other hand, argued that the incredible success of calculus and analytic methods in the 19th century had led mathematicians to focus on analytic functions and to ignore an important class of functions, namely, the discontinuous ones (Charles Hermite (1822-1901) called these ``Monsters'').  While his observation did in a vague sense foretell the development of discrete mathematics in the 20th century, Bugaev's specific mathematical direction here was ``arithmologie,'' a subject he devoted many papers to but which faded into oblivion when his students realized it was somewhat of a mathematical dead-end.    His students did, however, find another framework for carrying out Bugaev's vision of discontinuous functions, which they based on measure theory, and this was a driving force in the Moscow school for decades.

Bugaev's emphasis during his ICM talk was almost entirely philosophical rather than mathematical.  Discontinuous functions, he argued romantically, are morally strengthening because they free humankind from fatalism; discontinuity is a ``manifestation of independent individuality and autonomy.'' \cite[p. 68]{Naming}.  As Demidov deftly elaborates, analytic functions are completely determined by being prescribed in an arbitrarily small neighborhood of any point in their domain, and this in Bugaev's eyes is a mathematical expression of the idea of total determinism \cite{Demidov}.  In Svetlikova's words, Pavel Nekrasov (1853-1924) continued his mentor Bugaev's ideological pursuit and ``established a parallel between right-wing demands that the Russian Tsar be free from constitutional laws and the `arithmological' assumption that the power of natural laws was not absolute.'' \cite[p. 162]{Svetlikova}.  Whatever one makes of this, we see that Bely grew up in the environs of an influential, world-renowned, and politically reactionary (pro-Tsarist) math professor father whose mathematical work waded deeply into murky philosophical waters.

Let us turn now to Bely's novel, \emph{Petersburg}, which takes place in the eponymous Russian capital in 1905---the year and location, historically, of a failed political revolution that serves as the backdrop to the novel.  The primary protagonists are father and son.  Apollon Apollonovich Ableukhov is a conservative senator with an almost absurdist obsession with geometric ruminations.  It is widely recognized that much of his portrayal is intended as a caricature of Bugaev, and we shall discuss examples of this below.  Apollon's son, Nikolai, becomes entangled in a revolutionary plot and is given a bomb (referred to condescendingly as a sardine tin) that he is supposed to use on Apollon in a patricidal act of political terrorism.  

Since our goal in this paper is primarily to highlight and unwind Bely's mathematical prose, we shall not dwell on the plot, cast of characters, or literary significance of the novel (an interested reader will have no trouble finding endless writing on these matters!) and will instead only discuss what we need from the book as it comes up.  Before launching into our mathematical analysis, however, allow us an important point of context.  The Soviet literary scholar Konstantin Mochulsky, acclaimed for his landmark work on Dostoevsky, had this to say:
\begin{quote}
Andrei Bely's novel, \emph{Petersburg}, is the most powerful and artistically expressive of all his works.  It is a rendition of delirium unprecedented in literature; by means of subtle and complex literary devices, the author creates a separate world --- unbelievable, fantastic, monstrous, a world of nightmare and horror, a world of distorted perspectives, of disembodied people and living corpses.  [...] In order to understand the laws of this world, the reader must abandon at the threshold his own logical preconceptions.  Here common sense is abolished and causation enfeebled; here human consciousness is torn to shreds and explodes, like the hellish machine ``in the shape of a sardine tin.''  The ``topography'' of this world is strictly symmetrical: Petersburg with the straight lines of its prospects and the flatness of its squares, is perceived by the author as a system of ``pyramids,'' triangles, parallelepipeds, cubes and trapezoids.  This geometric space is peopled by abstract figures.  They move mechanically, act like automatons, and seem to read their lines from a script.  \cite[pp. 147-148]{Mochulsky}
\end{quote}

Some subsequent readings found more comedic satire in the oppressively dark imagery that pervades the novel than the unrelentingly nightmarish interpretation described above, though both responses remain common.  In describing Andrei Bely's personality, the influential Soviet writer Ilya Ehrenburg (1891-1967) articulated a connection between the author's mathematical obsessions and the peculiar emotional response his writing elicits in many readers, including Mochulsky \cite[p. 197]{Mochulsky}: ``Such a strange contradiction: the wild, fiery thought, while in his heart, instead of a burning coal, ice... Love and hate can carry people along with them, but not the madness of numbers, not the mathematics of the cosmos.  Bely's visions are full of magnificence and coldness.''  While Bely admittedly had a difficult and fractured relationship with his father, and this may have colored his attitude toward mathematics, we hope the reader takes away from our discussion in the following sections of this paper that Bely might also have found beauty in mathematics.  Whether he intended with his frequent ``Symbolist'' use of mathematics in \emph{Petersburg} to convey only the coldness described above or also the playful warmth those of us in the profession often feel when thinking about math is a matter left open for debate.

Note: All quotes and page numbers below refer to the 1978 Maguire and Malmstad English translation of \emph{Petersburg} \cite{Petersburg}.  See \cite{Translations} for a discussion of the different translations available, which is complicated by the fact that Bely published the novel serially in 1913-14, then in book form in 1916, then in a significantly edited and abbreviated form in 1922.  Maguire-Malmstad use the 1922 text, and they lay out their case for doing so in their introduction.

\section{Discontinuity}

Since discontinuity was such an important concept for Bugaev, it is natural to look for it in \emph{Petersburg}.  One passage that could be construed as a subtle wink from the author to his father is the following:
\begin{quote}
The hideous thing had been going on for twenty-four hours, or eighty thousand seconds, points in time.  Each moment was advancing.  And advancing on him was a moment which was somehow briskly spreading in circles, which was slowly turning into a swelling sphere.  This sphere was bursting, and his heel was slipping into voids.  Thus a wanderer in time was collapsing into the unknown and until ... a new instant. \cite[p. 217]{Petersburg}
\end{quote}

The imagery of expanding spheres is pervasive in the novel and has been discussed by various critics.  For a mathematical interpretation in terms of the Banach-Tarski paradox, see our companion paper \cite{Noah2}.  Our intent in quoting this passage here is instead to draw attention to Bely's description of time: he has discretized time into disconnected points/instants/moments.  His numerical conversion from hours to seconds psychologically amplifies this effect.

It is useful to recall that in 1900, Max Planck (1858-1947) suggested that radiation is quantized, meaning it comes in discrete amounts, then in 1913---the same year \emph{Petersburg} began to appear---Niels Bohr (1885-1962) constructed a theory of atomic structure based on quantum ideas, incorporating Plank's 1900 quantum hypothesis.  So the notion of the seemingly continuous universe secretly harboring discrete, discontinuous phenomena was quickly spreading.  Florensky, the student of Bugaev and classmate of Bely, carried his mentor's torch by arguing that intellectually the 19th century was a disaster, and the reason for this is that mathematics, specifically differential calculus, had created a monolith around the idea of continuity, leading society to believe that progress from one point to another without passing through all the intermediate points was impossible \cite[p. 87]{Naming}.  The above passage in \emph{Petersburg} allows time to move discretely through isolated temporal points, but Bely places this behavior in a rather delirious, nightmarish context; perhaps Bely's wink to his father and classmate was a sardonic one.

A more substantial manifestation of discontinuity in \emph{Petersburg} comes in the narrative structure of the book itself.  Bely deliberately rearranged sections of the book and blocks of text of various sizes to break the traditional linearity of prose.  In fact, he evidently confessed to using an elaborate ``cut and paste'' technique \cite[p. 101]{Alexandrov}.  Bely first applied this process in his 1902 work titled the ``Second Dramatic Symphony,'' and Mochulsky eloquently describes the result as follows \cite[p. 33]{Mochulsky}: ``the author tears apart the traditional structure of the narrative, the sequential flow of events joined by a causal and temporal link.  His prose is broken into pieces, smashed into fragments.''  Thus, Bely systematically introduces discontinuity into his prose, and this is one of the crucial modernist pillars underlying the success of \emph{Petersburg}.  In this respect, Bely is far from mocking his father's preoccupation with discontinuous functions---whether consciously or not, he has embraced his father's mathematical idea and innovatively transported it into his own literary creative process.

\section{Numbers as symbols}

Sprinkled throughout the text of \emph{Petersburg} are various passages that involve numbers in a way that is atypical in conventional literature.  At times Bely seems to almost meditate on the psychological impact of numbers or on the significance of counting---and countability.  This tends to set the mood of the novel and color our perceptions of various scenes; one could consider it to be a Symbolist use of mathematics (see also \cite[Chapter 5]{Columbia} for a discussion of the astral numerology in \emph{Petersburg}).

Right at the beginning of the novel, as Bely is first introducing us to the city of Petersburg (which many scholars have argued is the true protagonist of \emph{Petersburg}), the narrator states \cite[p. 2]{Petersburg}: ``It [Nevsky Prospect, the main street in the city] is delimited by numbered houses.  The numeration proceeds house by house, which considerably facilitates the finding of the house one needs.'' Although mild in terms of overt mathematical appearance, the early placement of this slightly odd remark suggests it may have a deceptively significant impact.  One interpretation is that Bely is quickly entering us into a land where mathematics underlies essentially everything.  Nekrasov, one of Bugaev's students and a known influence on Bely, maintained Bugaev's belief in the universal applicability of mathematics; one of Nekrasov's favorite quotations indicating this belief was the following \cite[p. 53]{Svetlikova}: ``Thou hast arranged all things by measure and number, and weight.''  This bears a certain resemblance to Bely's description of Nevsky Prospect.

While Florensky seems mostly to be emphasizing in his quote the idea that everything is quantifiable, Bely may actually have been reaching into a much deeper mathematical realm: countability, in the Cantorian sense.  As we learn in a modern undergraduate mathematics education, when faced with an infinite number of objects, being able to count, or enumerate, the objects constrains the size of the infinity.  Just ten pages after the preceding quote, Bely goes on to further describe the mathematical universe he is creating and in particular establishes the infinitude of his streets\footnote{The translators use the word ``prospect'' for street to keep closer to the Russian ``prospekt,'' since sound plays a tremendous role in all of Bely's writings.}:
\begin{quote}
The wet, slippery prospect was intersected by another wet prospect at a ninety-degree right angle.  At the point of intersection stood a policeman. [...] But parallel with the rushing prospect was another rushing prospect with the same row of boxes, with the same numeration, with the same clouds.  There is an infinity of rushing prospects with an infinity of rushing, intersecting shadows.  All of Petersburg is an infinity of the prospect raised to the nth degree.  Beyond Petersburg, there is nothing. \cite[pp. 11-12]{Petersburg}
\end{quote}
It appears that Bely is creating an infinite grid of orthogonal streets, so that the enumerated houses (the ``boxes'' he mentions) form a cartesian product of two countably infinite sets, and hence by Cantor's famous zig-zag argument the set of all houses in Petersburg is countably infinite.  Even Bely's ``nth'' power remains countable.  On the other hand, since there is nothing beyond Petersburg, his fictionalized city takes on an ``immeasurable immensity,'' as he later describes it.  Thus the mathematical tension uncovered by Cantor just a few decades earlier between the cardinalities of infinity is a background tension in Bely's novel.  Indeed, the $\aleph_0$ streets and houses give the inhabitants the local experience of an ordinary, navigable city, while the omniscient narrator makes the reader aware of the immeasurable vastness, and hence uncountable infinitude, of Petersburg.

That Bely was interested in the idea of infinity and countability comes as no surprise when one looks at the mathematical landscape in Russia at the time.  The book \emph{Naming Infinity} that we cited earlier \cite{Naming} is focused almost entirely on this story, so we shall only briefly recall it here.  The end of the 19th century produced what Hermann Weyl would later call the ``foundational crisis in mathematics.''  Cantor had rewritten the foundations of mathematics by introducing set theory, but in doing so he required mathematicians to accept infinity as an indispensable underpinning.  This shocked the mathematical world and directly contradicted a prevailing view by some that mathematics fundamentally describes the world and accordingly should be founded on the finite.  For many at the time, infinity was viewed as a mystical object, not a mathematical one.  The Axiom of Choice was formulated in 1904 by Ernst Zermelo (1871-1953) and led to intense debate.  The French school at the time (including Baire, Lebesgue, Borel, etc.) was strongly atheistic and actively steered mathematics away from questions of philosophy.  The Russian school in Petersburg (including Chebyshev, Lyapunov, Markov, etc.) maintained a similarly pragmatic outlook.

In direct contrast, the Moscow school, led by Bugaev, enthusiastically embraced the mystical implications of Cantor's work.  For instance, Florensky saw a close parallel between the religious act of naming God in prayer (particularly in the \emph{imiaslavie} movement) and Cantor's act of naming infinity in terms of cardinality.   Set theory came to underlie all of mathematics, and the Moscow school viewed mathematics as underlying the universe and all realms of human knowledge, so in their eyes this confluence of math and religion was inexorable.   The authors of \emph{Naming Infinity} argue that it was the Moscow school's embrace of mysticism in mathematics that allowed them to confront the foundational issues surrounding set theory at the time, and in doing so they launched down a very different mathematical path than the otherwise predominant French and Petersburg schools---and this schism played a large role in the ultimate growth of the Moscow school into one of the most active and influential in modern mathematical history.  Thus, while Bely was writing his novel about the attempted Russian revolution in 1905, his father's students were driving a mathematical revolution in Russia that literally took on biblical proportions.

The deep significance of ``naming'' that Florensky saw as a bridge between math and religion appeared quite explicitly in Bely's writing on the theoretical foundations of Symbolism.  In Alexandrov's words \cite[p. 103]{Alexandrov}, Bely argued that, ``naming a thing with a word causes that thing to come into existence.  All knowledge stems from naming, and knowing is impossible without words: before one can know something, that thing has to be defined, or, in effect, created with words.''  But we need not go as far as Cantor's infinite cardinalities to see Bely's fascination with naming numbers: even finite quantities can take on layers of Symbolist meaning.  For example, in one scene in \emph{Petersburg} \cite[p. 155]{Petersburg}, Nikolai describes his father's physical stature in terms of length and diameter.  By deliberately naming ``length'' instead of the more common and ostensibly equivalent ``height,'' Bely manages to reduce the father to something less than human---and this is taken even further with the use of ``diameter,'' which connects the father to the sphere symbol that is ubiquitous in the book.  In the following passage Bely insightfully plays with the import of a number's base 10 representation, which needless to say is just another form of it's \emph{name}:
\begin{quote}
Or--delirium that can be measured in digits.  Thirty zeroes---now that's a horror.  Yes, but cross out the numeral one and the thirty zeroes collapse.  There will be zero.  There is no horror in the numeral one.  The numeral one is a nonentity, a one and nothing more!  But one plus thirty zeros will make the hideous monstrosity of a quintillion.  A quintillion---oh, oh, oh!---dangles from a frail little stick.  The one of a quintillion repeats itself more than a billion billion times, themselves repeated more than a billion times.  Yes---Nikolai lived as a human numeral one, that is, as an emaciated little stick, running his course through time--- \cite[p. 226]{Petersburg}
\end{quote}

There are several other curious number-themed passages in \emph{Petersburg}.  At one point \cite[p. 22]{Petersburg}, the father Apollon Apollonovich ``had the sensation that his head was six times larger than it should be, and twelve times heavier than it should be.''  This odd phenomenon, and the peculiar numerical specificity in the description of the enlargement, sounds like a line from \emph{Alice's Adventures in Wonderland}!  (For a mathematical exegesis of that book, contextualized by the work of the Oxford math lecturer Charles Dodgson (1832-1898) who wrote it under the pseudonym Lewis Carroll, cf. \cite{AliceSolved}.)  Later on \cite[p. 26]{Petersburg}, the son Nikolai ``kept trying to catch at idle thoughts: thoughts of the number of books that would fit on the shelf of a bookcase [...].''  This passage makes one think of the enigmatic 1941 short story ``The Library of Babel'' by Jorge Luis Borges (1899-1986), and specifically the following sentence: ``Each wall of each hexagon [room in Borges' fictitious library] is furnished with five bookshelves; each bookshelf holds thirty-two books identical in format; each book contains four hundred ten pages; each page, forty lines; each line, approximately eighty black letters.''  (For a mathematical exegesis of \emph{this} story, cf. \cite{Borges}.)

Arguably the two most important number symbols for Bely are infinity and zero, which mathematically represent universal limits and visually are united by the circle motif in their typography.  We have seen infinity in Bely's geometric description of the streets of Petersburg, and here we find a mystical contemplation of ``zero'' in Nikolai's description of a fleeting vision:
\begin{quote}
I was growing, you see, into an immeasurable expanse, all objects were growing along with me, the room and the spire of Peter and Paul.  There was simply no place left to grow.  And at the end, at the termination---there seemed to be another beginning there, which was most preposterous and weird, perhaps because I lack an organ to grasp its meaning.  In place of the sense organs there was a ``\emph{zero}.''  I was aware of something that wasn't even a zero, but a \emph{zero minus something}, say five, for example. \cite[p. 182]{Petersburg}
\end{quote}
Providing a precise interpretation of this tantalizing quote has proven to be beyond our capabilities.

\section{Caricature and satire}

In this section we take a break from the heavy lifting of Symbolism, philosophy, and math history; we look here at Bely's skillful satire of the mathematical mind in \emph{Petersburg}.  Some of the quotes explored here are quite humorous; we expect them to be particularly entertaining to a mathematical readership, as Bely is remarkably successful at capturing the amusing eccentricities of our species (the mathematician, that is).  Bely did not confine this particular activity, mathematical satire, to his written work.  As Mochulsky notes \cite[p. 126]{Mochulsky}: ``Vyacheslav Ivanov [a fellow Symbolist poet] saw in Bely a resemblance to Gogol --- and laughed merrily when Bely, standing on the rug, would tell about his childhood, his father, the professors, imitating the persons of the Moscow eccentrics and performing parody scenes.''

Most of the mathematical satire in \emph{Petersburg} centers on the fictitious father Apollon Apollonovich and can with confidence be taken to represent Bely's real-life father, Bugaev.  Even though Apollon is a polician rather than a mathematician, the narrator's description of Apollon and his thoughts leaves no doubt about this link.  On page 14 we read that ``Apollon Apollonovich had abandoned himself to his favorite contemplation, cubes,'' and on page 72: ``As they ascended [the staircase], his legs formed angles, which soothed his spirit: he loved symmetry.''  On page 158: ``In the middle of the desk was a textbook entitled Planimetry.  Before going to bed, Apollon Apollonovich very often used to leaf through this little volume, so as to quiet the restless life inside his head with the most blissful outlines of parallelepipeds, parallelograms, cones, and cubes.''  There are also domestic similarities between Apollon and Bugaev (inadequacies as husbands and fathers), but that is beyond our scope.

Certain stereotypes about math professors evidently existed in much the same form a hundred years ago as they do today.  On page 32 we read: ``The tales of his [Apollon's] absentmindedness were legion.''  And on page 122: ``Apollon Apollonovich could not stomach face-to-face conversations, which, naturally, entailed looking each other in the eye.''  This latter quote reminds one of the well-known joke about shoes and an extroverted mathematician (if you haven't heard it before, ask a colleague or teacher).

Next, a caricature that, while equally unrealistic and superficial from a literary perspective, is significantly less so from a mathematical one:
\begin{quote}
We have cast an eye over this habitation, guided by the general characteristics which the senator was wont to bestow on all objects.  Thus: ---having found himself on one rare occasion in the flowering bosom of nature, Apollon Apollonovich saw: the flowering bosom of nature.  For us this bosom would immediately break down into its characteristics: into violets, buttercups, pinks; the senator would again reduce the particulars to a unity.  We would say, of course: ``There's a buttercup!'' ``There's a nice little forget-me-not....'' But Apollon Apollonovich would say simply and succinctly: ``A flower....''  Just between us: for some reason, Apollon Apollonovich considered all flowers the same, bluebells.  He would have characterized even his own house with laconic brevity, as consisting, for him, of walls (forming squares and cubes) into which windows were cut, of parquetry, of tables.  Beyond that were details. \cite[p. 21]{Petersburg}
\end{quote}
Here Bely has captured, albeit humorously, some of the genuine traits of mathematical thinking: abstraction, generalization, and categorization (or categorification!).  While the personification of these traits in Apollon is absurd to a comical extent, it is nonetheless an interesting depiction of mathematical cognition to find in a novel of any kind.

Throughout \emph{Petersburg}, Bely has merged the didactic pedantry and eccentric manners of the math professoriate with that of governmental bureaucracy.  This is exemplified most obviously with the dual father figure Apollon, but it reaches far beyond him.  At a party, on page 107: ``A little way off, the professor of statistics stumbled upon the zemstvo [local government] official, who was standing, bored, by the passageway.  He recognized him, smiled affably, and plucked at the button of his frock coat with two fingers as if grasping at his last means of salvation.  And now was heard: `According to statistical data ... the annual rate of consumption of salt by the average Dutchman ...'.''  The use of pronouns and passive voice here make it difficult to discern who said that last line to whom, but perhaps that is Bely's clever linguistic way of illustrating the similarity between these two (in his view dull and impersonal) professions.

In another scene, on page 249: ``He was still thinking about the same type of functionary.  People of that type always defend themselves with phrases like `as is well known,' when nothing is known, or `science teaches us,' when science does not teach.''  The translators' notes here mention that these were two favorite phrases of Joseph Stalin (1878-1953).  But the first one is also, even to this day, a common refrain among mathematicians hoping to sweep certain details under the rug---so whether intentional or not, Bely has struck another point of intersection between the quirks of mathematicians and unsavory politicians.  This satirical convergence of mathematics and politics reaches an even higher level of abstraction when it appears in the very language and thought process of the anonymous narrator.  A striking example of this occurs in the Prologue to the novel:
\begin{quote}
But if you continue to insist on the utterly preposterous legend about the existence of a Moscow population of a million-and-a-half, then you will have to admit that the capital is Moscow, for only capitals have a population of a million-and-a-half; but as for provincial cities, they do not, never have had, and never will have a population of a million-and-a-half.  And in conformance with this preposterous legend, it will be apparent that the capital is not Petersburg.  But if Petersburg is not the capital, then there is no Petersburg. It only appears to exist.  However that may be, Petersburg not only appears to us, but actually does appear---on maps: in the form of two small circles, one set inside the other, with a black dot in the center; and from precisely this mathematical point, which has no dimension, it proclaims forcefully that it exists: from here, from this very point surges and swarms the printed book; from this invisible point speeds the official circular. \cite[p. 2]{Petersburg} 
\end{quote}
This confounding, pseudo-logical verbal navigation could just as well be an improvised proof attempt at a blackboard in front of math students, or a dissembling, legalistic, political memo.

With a careful eye (or perhaps a creative imagination), one can find further play on mathematical proof-techniques embedded in the novel.  For instance, on page 130 Apollon Apollonovich mistakenly refers to the character Pavel Yakovlevich as Pavel Pavelovich.  After being introduced to a series of characters in the book whose given name matches their patronymic (Apollon Apollonovich, Konstantin Konstantinych, Hermann Hermannovich, Sergei Sergeyevich, Ivan Ivanovich), this could be interpreted as a sort of false-induction joke, like the young math student believing all odd numbers are prime since the statement holds for 3, 5, 7, ....   Another scene where the narrator mischievously indulges in logic-themed humor is the following:
\begin{quote}
The rooms were small.  Each was occupied by only one absolutely enormous object. In the minute bedroom the enormous object was the bed; in the minute bathroom, the bathtub; in the living room, the bluish alcove; in the dining room, the table and the sideboard; the enormous object in her husband's room was, it stands to reason, the husband. \cite[p 40]{Petersburg}
\end{quote}

\section{Recursion and fractals}

Bely has fun in his novel with the idea of recursion and self-similarity:
\begin{quote}
Apollon Apollonovich was like Zeus: out of his head flowed goddesses and genii.  One of these genii (the stranger with the small black mustache), arising as an image, had already begun to live and breath in the yellowish spaces.  And he maintained that he had emerged from there, not from the senatorial head.  This stranger turned out to have idle thoughts too.  And they also possessed the same qualities.  They would escape and take on substance.  And one fugitive thought was the thought that the stranger really existed.  The thought fled back into the senatorial brain.  The circle closed.  Apollon Apollonovich was like Zeus. \cite[pp. 20-21]{Petersburg}
\end{quote}
Not only does this passage implicitly bring up some of the wonderful paradoxical issues famously contemplated in Douglas Hofstadter's Pulitzer Prize-winning 1979 book \emph{G\"odel, Escher, Bach: An Eternal Golden Braid}, but it unites these with Bely's ubiquitous circle symbol (and in fact, the above passage itself forms a circle since it begins and ends with identical sentences).  For context, Russell's Paradox was discovered in 1901, not many years before Bely began \emph{Petersburg}, and this was just one of many bizarre phenomena that Cantor's set theory unleashed---and which, as we have discussed earlier, Bugaev's Moscow Mathematical School had embraced with religious fervor.  Several pages later, Bely takes his recursive loop one step further and even involves \emph{himself}, the author of \emph{Petersburg}:
\begin{quote}
In this chapter we have seen Senator Ableukhov.  We have also seen the idle thoughts of the senator in the form of the senator's house and in the form of the senator's son, who also carries his own idle thoughts in his head.  Finally, we have seen another idle shadow---the stranger.  This shadow arose by chance in the consciousness of Senator Ableukhov and acquired its ephemeral being there.  But the consciousness of Apollon Apollonovich is a shadowy consciousness because he too is the possessor of an ephemeral being and the fruit of the author's fantasy: unnecessary, idle cerebral play. \cite[p. 35]{Petersburg}
\end{quote}
From a stylistic perspective, this metafictional twist reaches both toward the future with postmodern literature and toward the past with Gogol's self-referencing narrator in his 1842 novel \emph{Dead Souls}.  From a mathematical perspective, one could say that it looks back at Cantor and forwards toward Kurt G\"odel (1906-1978), who based his 1931 proof of the famous Incompleteness Theorem on his so-called ``Self-Reference Lemma.''

Since recursion and self-similarity are closely related to fractals, let us conclude this section with the following passage, which is suffused with Bely's beloved revolutionary helical symbol:
\begin{quote}
The leaves stirred up from the spot.  They eddied in dry circles about the skirts of the greatcoat.  The circles narrowed and curled in ever more restless spirals.  The golden spiral whispered something and danced more briskly.  A vortex of leaves swirled, wound round and round, and moved off to the side somewhere, off to the side somewhere, without spinning.  \cite[p. 99]{Petersburg}
\end{quote}
Is Bely's ``golden spiral'' here a deliberate mathematical reference or a mere linguistic coincidence?  Perhaps from a Symbolist perspective the answer is irrelevant.

\section{Geometry and dimensions}

Some of the most striking mathematical writing in \emph{Petersburg} concerns geometry.  In addition to the ubiquity of spheres and cubes, which are mostly used symbolically, there are some truly elaborate, imaginative, geometric passages in which Bely lets his inner mathematician shine.  

Bely brings us into his geometric world early in the novel and sets the tone for what is to follow.  On pages 10-12, in the Chapter 1 section titled ``Squares, Parallelepipeds, Cubes'' we find the following vignettes (recall here that ``prospect'' means ``street''):
\begin{itemize}
\item ``He [Apollon] was cut off from the scum of the streets by four perpendicular walls.''
\item ``Proportionality and symmetry soothed the senator's nerves, which had been irritated both by the irregularity of his domestic life and by the futile rotation of our wheel of state.''
\item ``Most of all he loved the rectilineal prospect; this prospect reminded him of the flow of time between the two points of life.''
\item ``There the houses merged cubelike into a regular, five-story row.''
\item ``Inspiration took possession of the senator's soul whenever the lacquered cube cut along the line of the Nevsky: there the numeration of the houses was visible.''
\item ``At times, for hours on end, he would lapse into an unthinking contemplation of pyramids, triangles, parallelepipeds, cubes, and trapezoids.''
\end{itemize}

These little snippets on geometry culminate in the following flight of fancy on shape and space:
\begin{quote}
While gazing dreamily into that illimitability of mists, the statesman suddenly expanded out of the black cube of the carriage in all directions and soared above it.  And he wanted the carriage to fly forward, the prospects to fly to meet him---prospect after prospect, so that the entire spherical surface of the planet should be embraced, as in serpent coils, by blackish gray cubes of houses; so that all the earth, crushed by prospects, in its lineal cosmic flight should intersect, with its rectilineal principle, unembraceable infinity; so that the network of parallel prospects, intersected by a network prospects, should expand into the abysses of the universe in planes of squares and cubes [...]. \cite[p. 11]{Petersburg}
\end{quote}
Bely is playing here with two distinct but related tensions: that of linearity versus circularity (which we earlier saw him entwine in the spiral figure), and that of infinite versus finite (the earlier discussed topic of the foundational crisis in mathematics).  In the above passage Bely takes the linear streets of Petersburg and expands them to wrap around the surface of the spherical planet, which implicitly transforms them from lines into circles---even though the houses they delineate retain their cubical linearity during this transformation, naturally enough.  Moreover, science teaches us that planetary orbits are elliptical, yet Bely makes Earth's ``cosmic flight'' linear---as if at this very moment of Apollon's vision we turned off gravity and let our planet break free from its orbit, sailing in the direction of its tangent vector.   The spherical network of orthogonal parallel streets would then trace out linear paths in its galactic voyage.  This means the ``orbit'' of each house, and of each citizen of \emph{Petersburg}'s Petersburg, is a line in the three-dimensional cosmos.  Any two such lines are necessarily parallel, which means that in a Euclidean universe they never meet---but in projective space they meet at exactly one point at infinity: ``all the earth, crushed by prospects, in its lineal cosmic flight should intersect, with its rectilineal principle, unembraceable infinity.''  Did Bely deliberately intend to invoke projective geometry in his poetic mention of infinity here?  It is unclear, perhaps even unlikely, but arguably that only adds to the beauty and mystery of his writing, since this artistic mathematical vision he crafted is surprisingly plausible and consistent from a technical geometric perspective.

What of the tension between infinite and finite?  Bely earlier described the infinitude of his streets, and even here one can travel on any street unimpeded forever, yet the spherical geometry of the earth means that each street is of finite length.  This is quite similar to the situation of Borges' aforementioned Library of Babel.  In that story, inhabitants of the library can freely walk in any direction unimpeded forever, but by the narrator's description one infers that the library contains only finitely many books---which suggests that the library must exhibit some non-trivial global topology, the possibilities of which are informatively explored in William Bloch's book \cite{Borges}.   One option Bloch entertains is that Borges' library could be a non-orientable manifold, so that after traversing long enough in any direction, one returns to the same room but with all the books written in reverse.   Intriguing as it is, there is no particular evidence in Borges' writing that this is what was intended.  

On the other hand, this idea of reverse orientation does appear in \emph{Petersburg}!  On page 208, a spectral visitor declares: ``Yes, our spaces are not yours.  Everything there flows in reverse order.''  At this point the reader realizes the relation between two names that have been occurring throughout the book: Shishnarfne and Enfranshish (the translators do the best they can to preserve the palindromic symmetric when converting from Cyrillic).  This idea has precedent in the occult ``anthroposophy'' movement led by Rudolf Steiner (1861-1925) that Bely was absorbed in for some years.  As Alexandrov notes \cite[p. 119]{Alexandrov}: ``Steiner explains that elements of the spirit world, in addition to appearing as aspects of an individual's outer world, also appear as mirror images of what they really are: `When, for instance, a number is perceived, it must be reversed, as a picture in a mirror; 265 would mean here in reality 562'.''  While working on \emph{Petersburg} Bely had an intense falling out with Steiner, so whether he intends this other-worldly orientation reversal sincerely or satirically is not so clear.  Either way, the idea of \emph{Petersburg} existing in a non-orientable 3-manifold (or 4-manifold, as we soon discuss) is rather appealing to a mathematician.

When discussing satire, we noted that Bely conflates mathematics with conservative politics, most saliently in his fictionalized father.  One more example of this occurs early on: 
\begin{quote}
While dwelling in the center of the black, perfect, satin-lined cube, Apollon Apollonovich reveled at length in the quadrangular walls.  Apollon Apollonovich was born for solitary confinement.  Only his love for the plane geometry of the state had invested him in the polyhedrality of a responsible position.  \cite[p. 11]{Petersburg}
\end{quote}
This particular quote brings to mind the 1884 novella \emph{Flatland: A Romance of Many Dimensions}, by Edwin Abbott (1838-1926).  This story describes life in a 2-dimensional plane populated by line segments (women) and polygons (men); the more edges of a polygon, the higher the social status.  Equilateral triangles are lowly craftsmen; squares and pentagons are the ``gentlemen'' class including doctors and lawyers; hexagons are the lowest rank of nobility, and approximate circles comprise the high priest class.  Though we don't know whether Bely had read \emph{Flatland}, his senator's love for ``the plane geometry of the state'' and his ``polyhedrality of a responsible position'' is at the very least an amusing coincidence in light of Abbott's satirical Victorian world.  So, too, is Bely's assertion, ``After the line, the figure which soothed him more than all other symmetries was the square,'' which in \emph{Flatland} would humorously translate to: ``After women, the people who soothed Apollon Apollonovich the most were male doctors and lawyers.''  

In \emph{Flatland} the protagonist, a square, is visited by a sphere from the third dimension.  When passing through the planar Flatland, the sphere appears to Flatlanders first as a point (precisely when he is tangent to Flatland), then as a circle which enlarges then shrinks back to a point before disappearing.  The sphere is able to appear inside the homes of Flatlanders: no wall can keep him out, since he is not constrained to the plane of Flatland.  At one point in the story, the sphere enters the square's body and, after the former touches the inside of the latter's stomach, the square asserts that, ``a demoniacal laugh seemed to issue from within me.''   Then, as Banchoff observes \cite{Banchoff}, in Abbott's 1887 book \emph{The Kernel and the Husk} a discussion of spirits includes the following line: ``[...] a being of Four Dimensions, if such there were could come into our closed rooms without opening door or window, nay could even penetrate into, and inhabit our bodies.''  This perfectly extends the above visitation scene in \emph{Flatland} to one higher dimension, and an unmistakably similar idea occurs in \emph{Petersburg} involving Bely's demonic spirit character Shishnarfne:
\begin{quote}
A man of all three dimensions had entered the room.  He had leaned against the window and had become a contour (or, two-dimensional), had become a thin layer of soot of the sort you knock out of a lamp.  Now this black soot had suddenly smoldered away into an ash that gleamed in the moonlight, and the ash was flying away.  And there was no contour.  The whole material substance had turned into a phonic substance that was jabbering away.  But where?  It seemed to Alexander Ivanovich that the jabbering had now started up inside him.  ``Mr. Shishnarfne,'' said Alexander Ivanovich to space (for there was now no Shishnarfne).  And he chirruped as he answered himself: ``Petersburg is the fourth dimension which is not indicated on maps, which is indicated merely by a dot.  And this dot is the place where the plane of being is tangential to the surface of the sphere and the immense astral cosmos.  A dot which in the twinkling of an eye can produce for us an inhabitant of the fourth dimension, from whom not even a wall can protect us.''  \cite[p. 207]{Petersburg}
\end{quote}
Bely's grasp here of dimensionality and geometry is, like Abbott's, impressively accurate.  Again, whether the similarity is coincidence or inspiration we do not know.  A definite influence, however, comes once again from Rudolf Steiner; in Alexandrov's words \cite[p. 113]{Alexandrov}: ``The published version of one of Steiner's best-known lectures that he gave in 1911 has a sketch he drew showing the movement of `living beings on the Astral Plane' into a human head where `shadow-images' of these beings are reflected as thinking.''  A clever mathematical element Bely adds to both the Steiner and Abbott pictures is his use above of a tangent plane to relate a spherical world with a linear one.  And by placing Petersburg as the point of tangency, he sculpts it as both a mathematical and a spiritual point of contact between his two worlds.  Quite ingenious, really.

\section{Conclusion}
It is difficult, if not impossible, to name an important work of literature as heavily imbued with mathematics as Bely's \emph{Petersburg}.  This singular aspect of the novel has not escaped the attention of historian nor literary scholar, but nonetheless a ``close reading'' from a mathematician's eye seems not to have been undertaken previously.  In doing so, we find an amazingly rich array of mathematical manifestations and allusions.  Bely's father is caricatured for his pedantry, absent-mindedness, and penchant for abstraction, yet the structure of the novel itself reflects the father's universal faith in discontinuity.  Poetic shadows of Cantor's work on set theory, countability, and infinity appear in the novel and take on a Symbolist meaning in the context of the Moscow Mathematical School's religiously inspired and mystically driven work on set theory and measure theory.  We interpret Bely's fantastical description of the streets in Petersburg in terms of spherical and projective geometry.  We find striking similarities between Bely's treatment of a spiritual visitor from the fourth dimension and Abbott's famous \emph{Flatland}, in addition to a couple passages exhibiting a slight foreshadow to Borges' very mathematical short story ``Library of Babel.''

It is perhaps tempting to use these findings as clues that may help us crack Bely's code and decipher the true ``meaning'' of \emph{Petersburg}.  But this would be naive and injudicious.  Part of what distinguishes great modernist literature, particularly the four books Nabokov singled out, is their many levels and layers of meaning that allow critics to find nearly endless avenues of interpretation.  No one perspective is more correct than the others, but their union (so to speak) helps us understand the work in its proper context and reveals its literary depth and genius.  That one of the most important novels of the 20th century, written by the son of one of the most influential mathematicians of the late 19th century, is still revealing new secrets to us over a hundred years after its initial publication is astounding.  We hope that students of literature will appreciate our mathematical insight into Bely's work and world, and that students of mathematics will realize their subject reaches far beyond the classroom and textbooks in which they primarily encounter it.

\bibliographystyle{alpha}
\bibliography{bib}

\end{document}